\newcommand{\bull}{\vrule height .9ex width .8ex depth -.1ex}
 \newcommand{\ppp}{\hfill $\bull$ }

 \documentclass[11pt]{article}
 \author{ M.-L. Labbi}
 
   \title{On Two Natural Riemannian Metrics on a Tube}
   \date{}
   \usepackage{amsmath}
   \usepackage{amscd}

\begin{document}
   \maketitle
   \begin{abstract} 
During an operation of surgery on a Riemannian manifold and along a given 
embedded submanifold, (see \cite{{GroLaw},{Lab1},{Lab2}}), one needs to replace the (old)
metric induced by the exponential  map on a tubular neighborhood of the
 submanifold by  the Sasakian metric. So a good understanding of the
  behavior of these two metrics is important, this
is our main goal in this paper. In particular, we prove that
these two metrics are tangent up to the order one if and only if
 the submanifold is totally geodesic. In the case where the ambient space is an Euclidean space,
 we prove that the difference of these two metrics is quadratic  in the radius of the tube and 
depends only on the second fundamental form of the submanifold. Also the  case of spherical and hyperbolic space forms are studied.
   \end{abstract}
 \section{Statement of the results}
Let $(X,g)$ be a smooth Riemannian manifold of dimension $n+p$ and let 
$M$
be an embedded (compact) $n$-submanifold of $X$. Let
 $$ T_\epsilon=\{(x,v): x\in M,v\in N_xM \quad {\rm and}\quad 
g(v,v) <{\epsilon}^2\}$$
be a tube of radius $\epsilon$ around M, where $N_xM$ denotes the normal space
 to $M$
 at $x$.  It is well known that
 there exists $\epsilon_0>0$ such that the exponential map, 
${\rm exp}:T_\epsilon\rightarrow X$, is a
  diffeomorphism onto its image for all $\epsilon\leq \epsilon_0$. 
 We shall denote  by ${\rm exp}^*g$
  the pull back to $T_\epsilon$ of the metric $g$ on $X$.
\par
The normal sub-bundle $T_\epsilon$ can also be endowed with a second natural 
metric,
namely, the  Sasakian
 metric.
It is defined to be  the metric $h$ compatible with the normal connection 
of the normal (sub)bundle 
such that the natural projection $ \pi :(T_\epsilon,h)\rightarrow (M,g)$ is
 a Riemannian submersion.\par
  \medskip
In this paper we investigate the behavior of these two metrics  near the
 zero 
section of the normal bundle.
\par\noindent
  Let $(p,rn)$ be an arbitrary point in $T_\epsilon$, where $r<\epsilon$
 and $n$ is
 a unit normal
  vector to $M$ at $p$. 
 We shall  denote by  $A_n$    the shape operator 
  of the submanifold $M$ in the direction
 of  $n$.  
  \par\medskip\noindent
{\bf Theorem A.} {\sl Let $R$ denote the Riemann curvature $(0,4)$-tensor
 of $(X,g)$. Then  for $u_1,u_2\in T_{(p,rn)}T_\epsilon $, we have}
\begin{equation*}
\begin{split}
{\rm exp}^*g(u_1,u_2)&=h(u_1,u_2)-2g(A_n\pi_*u_1,\pi_*u_2)r\\
+&\bigl\{g(A_n\pi_*u_1,A_n\pi_*u_2)+R(\pi_*u_1,n,\pi_*u_2,n)+
{2\over 3}R(\pi_*u_1,n,Ku_2,n)\\
+&{2\over 3}R(\pi_*u_2,n,Ku_1,n)+
{1\over 3}R(Ku_1,n,Ku_2,n)\bigr\}r^2+O(r^3)
\end{split}
\end{equation*}
{\sl  In particular, }
 \begin{equation*}
\frac{d}{dr}_{\scriptstyle |_{r=0}}{\rm exp}^*g
 =\frac{d}{dr}_{\scriptstyle |_{r=0}}h-
2\pi^* (II_n). 
\end{equation*}
{\sl Where $II_n(u,v)=g(A_n\pi_*u,v)$ is the second fundamental 
form of $M$.}
 \par\medskip\noindent
{\sc Remark.} 
Note that in \cite{GroLaw}, at the beginning of the proof of Lemma 2 
page 430, it is claimed that
the  metrics ${\rm exp}^*g$ and $h$ are sufficiently close in the $C^2$-topology
as $r\rightarrow 0$. The same error is also in \cite{Lab1}. But this does  not affect the corresponding conclusions in both papers
(after minor changes), see \cite{Lab2}. \par\noindent
An alternative short way to notice this fact is as follows:\par\noindent
With respect to the metric  $h$, the zero section  $M \hookrightarrow T_{\epsilon}$ is totally geodesic (since for a Riemannian submersion the horizontal lift of a geodesic is a geodesic).
 But on the other side, the zero section  $M \hookrightarrow T_{\epsilon}$ is totally geodesic
for the metric ${\rm exp}^*g$ if and  only if $M$ is totally geodesic in $(X,g)$.\par\smallskip

 In the case when the ambient space $(X,g)$ is the Euclidean space 
${\bf R}^n$, we
  prove  the following simple formula relating the metrics
   ${\rm exp}^*g$ and $h$: \par\medskip\noindent
 {\bf Theorem B.} {\sl Let $M$ be an embedded submanifold in  the 
Euclidean
 space ${\bf R}^n$, then for $u_1,u_2\in T_{(p,rn)}T_\epsilon $, we have }
\begin{equation*}
{\rm exp}^*g(u_1,u_2)=h(u_1,u_2)-2g(A_n\pi_*u_1,\pi_*u_2)r
+g(A_n\pi_*u_1,A_n\pi_*u_2)r^2.
\end{equation*}
A similar result is proved for any space form, as follows:
\par\medskip\noindent
 {\bf Theorem C.} {\sl Let $M$ be an embedded submanifold in  a space form 
$(X,g)$ with curvature $k$, then for $u_1,u_2\in T_{(p,rn)}T_\epsilon $,
 we have }
\begin{equation*}
\begin{split}
{\rm exp}^*g(u_1,u_2)&=
\frac{\sin_k^2(r)}{r^2}h(u_1,u_2)-2\sin_k(r)\cos_k(r)g(A_n\pi_*u_1,u_2)\\
+&\sin_k^2(r)g(A_n\pi_*u_1,A_n\pi_*u_2)
+\{\cos_k^2(r)-
\frac{\sin_k^2(r)}{r^2}\}g(\pi_* u_1,\pi_* u_2)\\
+&
\{\frac{\sin_k(r)}{r}-1\}^2g(Ku_1,n)g(Ku_2,n).
\end{split}
\end{equation*}
\section{Preliminaries}
\subsection{The Sasakian Metric on the Normal Bundle}
Let $M$ and  $(X,g)$ be as above and let $\pi: \nu(M)\rightarrow M$ be the 
normal bundle of the embedding. Using the normal connection $\nabla$ of $\nu(M)$,
the tangent bundle $T(\nu(M))$ splits naturally to
$$T(\nu(M))={\cal V}\oplus {\cal H}.$$
Where ${\cal V}$ and ${\cal H}$ are respectively the vertical and
 horizontal bundles. Recall that at a given point $(p,v)\in \nu(M)$,
we have ${\cal V}_{(p,v)}=T_{(p,v)}\pi^{-1}(p)$, that is the tangent to
 the fiber over $p$. Hence,
using parallel displacement in the fiber, we can canonically 
identify the vertical space at $(p,v)$ with the fiber $\nu_p(M)$. 
Thus we get a map , called the 
connection map, 
$$K:T(\nu(M)) \rightarrow \nu(M)$$
It is the composition of the projection onto the vertical space followed by
a parallel displacement in the fiber as above. In particular, we have
$$K({\cal V}_{(p,v)})=\nu_p(M)\quad {\rm and}\quad K({\cal H}_{(p,v)})
=\{0\}.$$
More explicitly, if $u$ is a tangent vector at $t=0$  to a curve $(p(t),v(t))$  
 in $\nu(M)$, then
\begin{equation}\label{one}
Ku=\nabla_{\dot{p}}v(0).
\end{equation}
 On the other hand, a tangent vector $u$, as above, 
is horizontal if and only if $v(t)$ is $\nabla$-parallel along $p(t)$.
\par\noindent
 The Sasakian metric on $\nu(M)$ is  defined by
\begin{equation}\label{two}
h(u_1,u_2)=g(Ku_1,Ku_2)+g(\pi_*u_1,\pi_*u_2).\end{equation}
Note that clearly $\pi: (\nu(M),h)\rightarrow (M,g)$ is then a
Riemannian submersion.
\subsection{Jacobi Fields and the Exponential Map}
Let $r$ be positive and let $n$ be a unit normal vector at $p\in M$.
Let  $u\in T_{(p,rn)}\nu(M)$, 
then $u={d\over dt}_{|_{t=0}}
(p(t),rn(t))$. Consider $U(s)={d\over dt}_{|_{t=0}}(
p(t),sn(t))$. It is a vector field along the curve $c(s)=(p,sn)$ such that 
$U(r)=u$.
Next, set 
\begin{equation*}
Y(s)={\rm exp}_*U(s).
\end{equation*}
\par\medskip\noindent
{\sc Fact:} {\sl The vector field $Y(s)$ is a Jacobi field in $(X,g)$ 
such
that
\begin{equation}\label{three}\begin{split}
Y(o)=\dot{p}(0)=\pi_*(u)\in T_pM,\\
{D\over ds}Y(0)={1\over r}Ku-A_n(\pi_*(u)).
\end{split}\end{equation}
Where $D$ and $A_n$ denote respectively the Riemannian connection of $(X,g)$
and the shape operator of $M$.
}\\
\par\medskip\noindent
{\sc Proof:}
Let $\xi(s)={\rm exp}_{\nu}c(s)={\rm exp}_{\nu}(p,sn)$ be the unit speed
 geodesic in $(X,g)$ normal to $M$ with 
$\xi(0)=p\in M$ and $n=\dot{\xi}(0)\in \nu_p(M)$.
Remark that
 $$Y(s)={\rm exp}_*U(s)=
{d\over dt}_{|_{t=0}}{\rm exp}_\nu(p(t),sn(t)).$$
 The vector field $Y(s)$ is then  generated from a variation
of geodesics, in $(X,g)$,  normal to $M$, namely, 
$c(t,s)={\rm exp}_\nu(p(t),sn(t))$. Then
 $Y(s)$ is a Jacobi field  along $\xi(s)$. Furthermore, we have\\ 
\begin{equation*}
Y(o)=\dot{p}(0)=\pi_*(u)\in T_pM.\end{equation*}
 Also, using (\ref{one}) we get
\begin{equation*}
\begin{split}
{D\over ds}Y(0)=&{D\over ds}
\bigl({d\over dt}_{|_{t=0}}{\rm exp}_\nu(p(t),sn(t))\bigr)(0)
={D\over dt}\bigl({d\over ds}_{|_{s=0}}{\rm exp}_\nu(p(t),sn(t))\bigr)(0)\\
=&{D\over dt}n(t)(0)
={1\over r}{D^{\bot}\over dt}rn(t)(0)+{D^{T}\over dt}n(t)(0)\\
=&{1\over r}Ku-A_n(\pi_*(u)).
\end{split}
\end{equation*}
This completes the proof.\ppp
\par\noindent
Finally, note the following remarks:
\begin{equation*}\begin{split}
\pi_*(U(s))=Y(0)=\pi_*(u)=\dot{p}(0),\\
K(U(s))=\nabla_{\dot{p}}sn(t)={s\over r}\nabla_{\dot{p}}rn(t)=
{s\over r}Ku.
\end{split}\end{equation*}
\section{Proof of Theorems}
Let $u_1,u_2\in T_{(p,rn)}\nu(M)$, then 
\begin{equation*}
{\rm exp}^*g(u_1,u_2)=g({\rm exp}_*u_1,{\rm exp}_*u_2)=
g({\rm exp}_*U_1(r),{\rm exp}_*U_2(r))=g(Y_1(r),Y_2(r)),
\end{equation*}
where, for $j=1,2$, $Y_j(s)={\rm exp}_*U_j(s)$ and $U_j$ is  the vector field
 associated to the vector
$u_j$ as above.\\
\subsection{Proof of Theorem A}
It results from the Jacobi equation that for every $k\geq 0$ and for 
$j=1$ or $2$,
 we have
\begin{equation*}
D_{\dot{\xi}}^{k+2}Y_j(o)=-\sum_{i=0}^kC_i^k(D_{\dot{\xi}}^{k-i}R)
(D_{\dot{\xi}}^{i}Y_j(0),n)n
\end{equation*}
In particular, using  (\ref{three}), we have
\begin{equation*}
\begin{split}
D_{\dot{\xi}}^{2}Y_j(o)&=-R(\pi_*u_j,n)n,\\
D_{\dot{\xi}}^{3}Y_j(o)&=-D_{\dot{\xi}}R(\pi_*u_j,n)n-{1\over r}
R(Ku_j,n)n+R(\pi_*u_j,n)n\\
D_{\dot{\xi}}^{4}Y_j(o)&=-D_{\dot{\xi}}^2R(\pi_*u_j,n)n-{2\over r}
D_{\dot{\xi}}R(\pi_*u_j,n)n+2D_{\dot{\xi}}R(A_n\pi_*u_j,n)n\\
&+
R(R(\pi_*u_j,n)n,n)n
\end{split}
\end{equation*}
Next, the Taylor expansion of $g(Y_1(r),Y_2(r))$ shows that
\begin{equation*}
\begin{split}
{\rm exp}^*g(u_1,u_2)&=g(Y_1(r),Y_2(r))\\
&=g(\pi_*u_1,\pi_*u_2)-2g(A_n\pi_*u_1,\pi_*u_2)r\\
&+\bigl\{2R(\pi_*u_1,n,
\pi_*u_2,n)
+{2\over {r^2}}g(Ku_1,Ku_2)
+2g(A_n\pi_*u_1,A_n\pi_*u_2)
\bigr\}{r^2\over 2!}\\
&+
\bigl\{{4\over r}R(\pi_*u_1,n,Ku_2,n)
+{4\over r}R(Ku_1,n,\pi_*u_2,n)+O(1)
\bigr\}{r^3\over 3!}\\
+&\bigl\{{8\over {r^2}}R(Ku_1,n,Ku_2,n)+O({1\over r})\bigr\}{r^4\over 4!}+...
\end{split}
\end{equation*}
Consequently, using (\ref{two}) we get 
\begin{equation*}
\begin{split}
{\rm exp}^*g(u_1,u_2)&=h(u_1,u_2)-2g(A_n\pi_*u_1,\pi_*u_2)r
+\bigl\{g(A_n\pi_*u_1,A_n\pi_*u_2)\\
+& R(\pi_*u_1,n,\pi_*u_2,n)+
{2\over 3}R(\pi_*u_1,n,Ku_2,n)+{2\over 3}R(\pi_*u_2,n,Ku_1,n)\\
+&
{1\over 3}R(Ku_1,n,Ku_2,n)\bigr\}r^2+O(r^3).
\end{split}
\end{equation*}
This completes the proof. \ppp
\subsection{Proof of Theorems B and C}
Here, we suppose the manifold $(X,g)$ is with constant sectional curvature
$k$,
then the Jacobi equation for $Y_j$, $j=1$ or $2$, becomes
\begin{equation}\label{jaceq}
Y_j''(s)+kY_j(s)-kg\bigl(Y_j(s),\dot{\xi}(s)\bigr)\dot{\xi}(s)=0.
\end{equation}
Next, note that
$$Y_j(s)=g\bigl(Y_j(s),\dot{\xi}(s)\bigr)\dot{\xi}(s)+Y_j^{\bot}(s),$$
and it is easy to check that
$$g(Y_j(s),\dot{\xi}(s))=g(Y_j(0),\dot{\xi}(0))+g(Y'_j(0),\dot{\xi}(0))s.$$
On the other hand, we also have
$$(Y_j^{\bot})^{\prime\prime}(s)+R\bigl(Y_j^{\bot}(s),\dot{\xi}(s)\bigr)
\dot{\xi}(s)=0,$$
then the vector field $Y_j^{\bot}(s)$ satisfies
\begin{equation*}
\begin{split}
(Y_j^{\bot})''(s)&+kY_j^{\bot}(s)=0,\\
Y_j^{\bot}(0)&=Y_j(0),\\
(Y_j^{\bot})'(0)&=Y_j'(0)-g(Y_j'(0),n)n.
\end{split}
\end{equation*}
The solutions of this differential equation  are in terms of
parallel
translation $\tau_s$ along $\xi(s)$ as follows
\begin{equation}
Y_j^{\bot}(s)=\cos_k(s)\tau_s(Y_j^{\bot}(0))+
\sin_k(s)\tau_s\bigl((Y_j^{\bot})'(0)\bigr).
\end{equation}
Where, $\cos_k(s)={d\over ds}\sin_k(s)$ and
\begin{equation}
\sin_k(s)=
\begin{cases}
 {\sin\sqrt{k}s\over \sqrt{k}} &\text {if $k>0$}\\
s&\text{if $k=0$}\\
{\sinh\sqrt{|k|}s\over \sqrt{|k|}} &\text {if $k<0$}
\end{cases}
\end{equation}
 Consequently, after using formula (\ref{three}), the Jacobi fields 
$Y_j$ are explicitly given by
\begin{equation*}\begin{split}
Y_j(s)=\frac{s}{r}g(Ku_j,n)\dot{\xi}(s)+&\cos_k(s) \tau_s(\pi_*u_j)\\
+&
\sin_k(s)\tau_s\bigl\{\frac{1}{r}Ku_j-A_n\pi_*u_j-g(\frac{1}{r}Ku_j,n)n
\bigr\}
\end{split}
\end{equation*}
Finally, a direct computation shows that
\begin{equation*}
\begin{split}
{\rm exp}^*g(u_1,u_2)&=g(Y_1(r),Y_2(r))\\
=&\frac{\sin_k^2(r)}{r^2}h(u_1,u_2)-2\sin_k(r)\cos_k(r)g(A_n\pi_*u_1,u_2)\\
+&\sin_k^2(r)g(A_n\pi_*u_1,A_n\pi_*u_2)
+\{\cos_k^2(r)-
\frac{\sin_k^2(r)}{r^2}\}g(\pi_* u_1,\pi_* u_2)\\
+&
\{\frac{\sin_k(r)}{r}-1\}^2g(Ku_1,n)g(Ku_2,n).
\end{split}
\end{equation*}
This completes the proof. \ppp
\par
\par\noindent
M.-L. Labbi\\
Department of Mathematics,\\ 
College of Science,\\
 University of Bahrain, \\
32038 Isa Town,\\
Kingdom of Bahrain.\\
E-mail: labbi@sci.uob.bh 
\end{document}